%% file: paper.tex
\documentclass{article}

\usepackage{float}
\usepackage{mathtools}
\usepackage{cases}

\usepackage{graphicx}
\usepackage{amsmath,amsthm,bm,mathrsfs}
\usepackage[overload]{empheq}
\usepackage{xcolor}

\input standard.tex

\numberwithin{equation}{section}

\makeatletter
\let\@fnsymbol\@arabic
\makeatother

\begin{document}

\title{Asymptotic Analysis of a Two-Phase Model of Solid Tumour Growth}
\author{{\sc Andrea Genovese de Oliveira}\thanks{Departamento de Matem\'atica, Universidade de Bras\'ilia, Brazil, email: {\tt andreagenovese@unb.br}.}
\quad
{\sc John R. King}\thanks{School of Mathematical Sciences, University of Nottingham, United Kingdom, email: {\tt John.King@nottingham.ac.uk}.}}

\pagestyle{headings}
\markboth{A. GENOVESE DE OLIVEIRA}{\rm ASYMPTOTIC ANALYSIS OF TUMOUR MODEL}
\maketitle


\begin{abstract}
{\noindent We investigate avascular tumour growth as a two-phase process consisting of cells and liquid.  Based on the one-dimensional continuum moving-boundary model formulated by (\textit{Byrne, King, McElwain, Preziosi, Applied Mathematics Letters, 2003, 16, 567-573}), we defined boundary conditions for the analogous model of tumour growth in two dimensions.  We investigate linear stability of one dimensional time-dependent solution profiles in the moving-boundary formulation of a limit case (with negligible nutrient consumption and cell drag).  For this, we obtain an asymptotic limit of the two-dimensional perturbations for large time (in the case where the tumour is growing) by using the method of matched asymptotic approximations.  Having characterised an asymptotic limit of the perturbations, we compare it to the time-dependent solution profile in order to analytically obtain a condition for instability.  Numerical simulations are mentioned.  

\vspace{.3cm}

}
\end{abstract}


\section{Introduction}

We explore the two-phase model of avascular tumour growth proposed by Byrne, King, McElwain and Preziosi in 2003 \cite{byrne2003two}.  This is one of the seminal multiphase models of tumour growth.  By considering two limit cases of this model, the authors were able to draw parallels to previous models of tumour growth such as \cite{adam1986simplified}, \cite{gatenby1996reaction}, \cite{greenspan1972models}, and \cite{sherratt1992oncogenes}, .  On the other hand, it has been extended to three and four phase models of vascular tumour growth in \cite{breward2003multiphase} and \cite{hubbard2013multiphase}, respectively, using similar arguments as the ones used in this paper.  

Multiphase models consider the simultaneous movement of materials in different states (liquids, gases and solids) or with different chemical properties, such as viscosity and heterogeneity.  In tumour modelling, these phases can be, among others, healthy and cancer cells, extra-cellular material, and blood vessels, which will clearly have different properties.  Additionally, these phases may not be all included simultaneously in every model, depending on the complexity needed.  Indeed, this makes multiphase formulations conveniently constructive.  

We will approach this model through the study of stability, which is important in tumour growth models.  Stability can represent a non-malignant dormant tumour that stabilises at a certain size and can potentially be moved by redirecting its nutrient supply.  On the other hand, instability may represent a malignant tumour that grows uncontrollably, disintegrates or even breaks into pieces.  In the study of avascular tumours, particularly, the distinction between dormant and malignant tumours is a central point of research and this motivates the use of linear stability and perturbation methods to study the model proposed here.       

More models of tumour growth can be found in the following reviews: \cite{araujo2004history,roose2007mathematical}

\section{Model formulation}

The equations we analyse are those of the two phase model of tumour growth defined in \cite{byrne2003two}, namely 
\begin{alignat}{2}[left={\empheqlbrace}]
\dfrac{\partial \alpha}{\partial t}  & = S_c(\alpha,C)- \nabla \cdot (\mathbf{v}_c \alpha), \label{eq:2.1} 
\\
0 & = \nabla \cdot(\alpha\mathbf{v}_c+(1-\alpha)\mathbf{v}_w),\label{eq:2.2}
\\
\mathbf{0} & = \nabla \cdot(-p \mathbf{I} -\alpha \Sigma_c(\alpha)\mathbf{I} + \mu_c \alpha ( \nabla \mathbf{v}_c + \nabla \mathbf{v}_c^T) +  \lambda_c \alpha (\nabla \cdot \mathbf{v}_c)\mathbf{I}),\label{eq:2.3}
\\
\mathbf{0} & = (1-\alpha) \nabla p + k(\alpha)(\mathbf{v}_w- \mathbf{v}_c),\label{eq:2.4}
\\
0 & =\nabla^2 C - Q_c(\alpha, C),\label{eq:2.5}
\end{alignat}
%
wherein \eqref{eq:2.1} represents the conservation of mass equation for the cell phase (with cell volume fraction $\alpha$, cell velocity $\mathbf{v}_c$, net birth rate $S_c$), \eqref{eq:2.2} represents overall mass conservation (with water volume fraction $1-\alpha$, water velocity $\mathbf{v}_w$), \eqref{eq:2.3} is the overall momentum conservation equation (with water pressure $p$, cell extra pressure $\Sigma_c$, cell shear viscosity $\mu_c$, and cell bulk viscosity $\lambda_c $), \eqref{eq:2.4} is the water momentum equation (with interphase drag coefficient $k$) and \eqref{eq:2.5} is the (quasi-steady) nutrient consumption equation (with consumption rate $Q_c$).  We focus here on the two-dimentional case, with spatial coordinate $(x,y)$ and time $t$.

The following expressions were used for the functions $S_c(\alpha,C)$, $k(\alpha)$, and $Q_c(\alpha,C)$:
\begin{align*}
S_c(\alpha, C)  = \Big(\dfrac{s_0C}{1+s_1C}\Big)\alpha(1-\alpha) - \Big(\dfrac{s_2+s_3C}{1+s_4C}\Big)\alpha, \quad s_i>0,\quad  i=0\ldots 4 \\
k(\alpha)  = k_0 \alpha(1-\alpha), \quad k_0 >0, \qquad Q_c(\alpha, C)  = \dfrac{Q_0 C \alpha}{1+Q_1C}, \quad Q_0, Q_1 >0.
\end{align*}
The function $\Sigma_c(\alpha)$ is defined as the difference between the pressures of the two phases, chosen as in \cite{byrne2003two}:
$$
\Sigma_c(\alpha) = p_c - p_w =\dfrac{\hat{\Sigma}_c |\alpha - \alpha^\ast|^{r-1}}{(1-\alpha)^q}(\alpha-\alpha^\ast) H(\alpha-\alpha_{min}),$$
where $H$ is the Heaviside function and $q,r,\hat\Sigma_c>0$.  Additionally, $0<\alpha_{min}<\alpha^\ast<1$ with $\alpha^\ast$ being the natural cell density. Here, the cell pressure $p_c$ represents the isotropic stresses associated with cell-cell interactions and the fluid pressure $p_w$ is the hydrodynamic pressure in the water.  In summary, the equation above states that the difference in pressure between the two phases may depend on the cell concentration.

Additionally, we define the outer boundary of the tumour at time $t$ by
\begin{equation*}
\Gamma(t) :=: \{(x,y)\in \mathcal{R}^2: x=R(y,t)\}
\end{equation*}
such that the kinematic condition reads
\begin{equation} \label{eq:2.6}
\dfrac{\partial \mathbf{R}}{\partial t}(y,t) \cdot \mathbf{n}(y,t) = \mathbf{v}_c(\mathbf{R}(y,t),t)\cdot \mathbf{n}(y,t),
\end{equation}
where $\mathbf{n}$ is the outward unit normal to the curve $\mathbf{R}(y,t)=(R(y,t),y)$.  At the moving boundary $\Gamma(t)$ we also impose the following boundary conditions:
\begin{equation} \label{eq:2.7}
[- \Sigma_c(\alpha)\mathbf{I} + \mu_c  ( \nabla \mathbf{v}_c + \nabla \mathbf{v}_c^T) +  \lambda_c  (\nabla \cdot \mathbf{v}_c)\mathbf{I})]\cdot \mathbf{n} = \mathbf{0}, \qquad p=0, 
\qquad C=C_\infty.
\end{equation}
These represent zero stress in each phase and fixed nutrient concentration at the free surface.  Additionally, at the inner boundary $x=0$ we impose the following conditions:
\begin{equation} \label{eq:2.8}
\mathbf{v}_c=0, \qquad v_w^1 :=\mathbf{v}_w \cdot \mathbf{e}_1=0, \qquad \dfrac{\partial C}{\partial x} = 0,
\end{equation}
where $\mathbf{e}_1$ is the unit vector in the $x$ direction. This corresponds to assuming that $x=0$ is a rigid impermeable boundary.    

Notice that the boundary conditions imposed on the moving boundary can be reduced to the conditions required in the one dimensional model in \cite{byrne2003two} by taking 
$$\mathbf{R}(y,t)=(R(t),y), \qquad   \mathbf{n}=\mathbf{e}_1. 
$$
Additionally, recall that, in the one dimensional model, the outer boundary conditions on $v_c$ and $p$ were considered to be derived from $1$-D versions of 
$$\boldsymbol\sigma_c \cdot \mathbf{n} = 0, \qquad\boldsymbol\sigma_w \cdot \mathbf{n} = 0,$$
which are equivalent to the first two conditions defined on the two dimensional moving boundary.  We observe that \eqref{eq:2.4} and \eqref{eq:2.7} imply that 
\begin{equation} \label{eq:2.9}
(\mathbf{v}_w-\mathbf{v}_c) \cdot \boldsymbol{\tau}=0
\end{equation} 
on $\Gamma(t).$

To complete the necessary notation, let 
\begin{equation*}
\Omega(t) :=: \{(x,y) \in \mathcal{R}^2: 0<x<R(y,t)\}
\end{equation*}
denote our domain of interest for each $t>0$.  
  
\section{Limit case}\label{sec:3}

Motivated by results obtained in the simpler one dimensional case, we will analyse a specific limit case of this system in which nutrient consumption and cell drag are taken to be negligible.  This will allow us to obtain some analytically tractable simplifications of the full model.  

In the one dimensional formulation this limit case is obtained by taking $Q_c(\alpha,C) \equiv 0$ and $k(\alpha)\equiv 0$ in the original formulation.  The former raises no difficulties, simply implying that $C \equiv C_\infty$, (i.e. nutrient-rich behaviour) throughout $\Omega(t)$, but the latter requires more care.   If $k(\alpha) \equiv 0$ in the two dimensional formulation, then \eqref{eq:2.4} and \eqref{eq:2.7} imply 
$p\equiv 0$, only one equation, from the vector system \eqref{eq:2.4}.  To avoid this loss of information, we must reinstate the small $k(\alpha)$ taking $k(\alpha) = k_0 \alpha(1-\alpha)$ as in \cite{byrne2003two}, with the constant $k_0>0$ small.  Then,
\begin{equation*}
0 = \nabla p + k_0 \alpha  (\mathbf{v}_w-\mathbf{v}_c).
\end{equation*}
and hence
\begin{equation*}
0 = \nabla \times (\alpha( \mathbf{v}_w-\mathbf{v}_c))
\end{equation*}
and thus we have obtained the second equation (and extra two equations in three dimensions) required.

Therefore,  $p\equiv 0$ and $(\alpha, C, \mathbf{v}_c, \mathbf{v}_w, \mathbf{R})$ satisfy
\begin{alignat}{2}[left={\empheqlbrace}]
\dfrac{\partial \alpha}{\partial t} & = S_c(\alpha,C) - \nabla \cdot (\mathbf{v}_c \alpha),
\\
0 & = \nabla \cdot(\alpha\mathbf{v}_c+(1-\alpha)\mathbf{v}_w),
\\
0 & = \nabla \cdot(-\alpha \Sigma_c(\alpha)\mathbf{I} + \mu_c \alpha ( \nabla \mathbf{v}_c + \nabla \mathbf{v}_c^T) +  \lambda_c \alpha (\nabla \cdot \mathbf{v}_c)\mathbf{I}),
\\
0 &= \nabla \times [\alpha( \mathbf{v}_w-\mathbf{v}_c)]
\\
0 & = \nabla^2 C, \qquad\qquad\qquad\qquad (x,y)\in \Omega(t), \quad t\in(0,\infty), 
\end{alignat}
with kinematic condition \eqref{eq:2.6} and other boundary conditions as in \eqref{eq:2.7} and \eqref{eq:2.8}.  

Motivated by results obtained in the one dimensional case, described in \cite{genovese2017asymptotic}, we consider the stability of the following non-homogeneous solution of the limit case defined in the previous section.  Let the constants $\alpha=\alpha_h\in (0,1)$ and $\lambda_2$ be defined by 
\begin{equation*}
\lambda_2 = \dfrac{\Sigma_c(\alpha_h)}{\hat\mu_c} =  \dfrac{ S_c(\alpha_h, C_\infty)}{\alpha_h}
\end{equation*}
and 
$$\mathbf{v}_{c}(x) = \lambda_2 x \hspace{.05cm} \mathbf{e_1}\qquad 
\mathbf{v}_{w}(x) = -\dfrac{\lambda_2 \alpha_h}{1-\alpha_h}x \hspace{.05cm} \mathbf{e_1}, \qquad C=C_\infty.$$
We can also define the two-dimensional moving boundary $\Gamma_\ast (t)$ and domain $\Omega_\ast(t)$, via 
\begin{equation*}
\mathbf{R}_\ast(t,y) = (R_\ast(t),y) = (R_0 e^{\lambda_2 t},y).
\end{equation*}

\section{Linear stability analysis} 

We shall consider two dimensional perturbations on this one dimensional solution.  For this purpose, we linearise about this one dimensional solution and obtain a system for the respective two-dimensional perturbations
$$
\begin{array}{c}
\hat\alpha(x,y,t), \qquad \hat C(x,y,t), \qquad \hat p(x,y,t),\qquad \mathbf{\hat R}(y,t) =(\hat R(y,t),0)\\
\mathbf{\hat{v}}_{c}(x,y,t)=(\hat v^1_{c}(x,y,t),\hat v^2_{c}(x,y,t)) ,\quad \mathbf{\hat{v}}_{w}(x,y,t)=(\hat v^1_{w}(x,y,t),\hat v^2_{w}(x,y,t)).
\end{array}
$$
Then, we fix the boundary.  For this, we make the following change of variable:
\begin{equation*}
x =  R_\ast(t)\xi
\end{equation*}
and therefore
\begin{align*}
\dfrac{\partial}{\partial t} & \longrightarrow  \dfrac{\partial}{\partial t} - \dfrac{\xi}{R_\ast(t)} \dfrac{dR_\ast}{dt} \dfrac{\partial}{\partial \xi}
\\[.1cm]
\dfrac{\partial}{\partial x} & \longrightarrow  \dfrac{1}{R_\ast(t)} \dfrac{\partial}{\partial \xi}
\\[.1cm]
\dfrac{\partial}{\partial y} & \longrightarrow  \dfrac{\partial}{\partial y}.
\end{align*}

Having done this, we noticed that it is not possible to separate variables in $\xi$ or $t$ as both $\xi$ and $t$ are present in the coefficients in ways that they cannot be extracted by simple exponentiation.  However, it is possible to separate variables in $y$.  

Indeed, consider now
\begin{equation*}
\hat\alpha(\xi,y,t)= e^{i\kappa y}\tilde{\alpha}(\xi,t)+ c.c.
\end{equation*}
\begin{equation*}
\hat C(\xi,y,t)= e^{i\kappa y}\tilde{C}(\xi,t)+ c.c.\qquad\hat p(\xi,y,t)= e^{i\kappa y}\tilde{p}(\xi,t)+ c.c.,
\end{equation*}
\begin{equation*}
\mathbf{\hat{v}}_{c}(\xi,y,t)= e^{i\kappa y}\mathbf{\tilde{v}}_c(\xi,t) + c.c.= e^{i\kappa y}(\tilde v_c^1(\xi,t),\tilde v_c^2(\xi,t))+ c.c.,
\end{equation*}
\begin{equation*}
\mathbf{\hat v}_{w}(\xi,y,t)= e^{i\kappa y}\mathbf{\tilde{v}}_w(\xi,t)+ c.c.= e^{i \kappa y}(\tilde v_w^1(\xi,t),\tilde v_w^2(\xi,t))+ c.c.
\end{equation*}
and
$$ \mathbf{\hat R}(y,t) = (\hat R(y,t), 0)+ c.c.=e^{i\kappa y} (\tilde R(t),0)+ c.c,$$
where $c.c.$ stands for the complex conjugate of its preceding term and $\kappa\in \mathcal{R}$ is arbitrary corresponding to the wavelength of the perturbation in question.    

This reduces the system to
\begin{align}\begin{cases}
\label{eq:limit2d}
\dfrac{\partial \tilde \alpha }{\partial t} &=  \Big[\dfrac{\partial S_c}{\partial \alpha}(\alpha_h, C_\infty) -\lambda_2 \Big]\tilde \alpha  - \dfrac{\alpha_h }{R_\ast} \dfrac{\partial \tilde v^1_{c}}{\partial \xi} -\alpha_h i\kappa \tilde v_c^2
\\[.3cm]
0 &= \dfrac{\lambda_2}{1-\alpha_h}\Big(1+\xi \dfrac{\partial }{\partial \xi}\Big) \tilde  \alpha  +  \dfrac{\alpha_h}{R_\ast} \dfrac{\partial  \tilde  v_c^1}{\partial \xi }  +\dfrac{ 1-\alpha_h}{R_\ast} \dfrac{\partial \tilde  v_w^1 }{\partial \xi }\\[.3cm]
& +\alpha_h i\kappa \tilde  v_c^2 + (1-\alpha_h) i\kappa  \tilde  v_w^2
\\[.3cm]
0 & = - \dfrac{\Sigma_c'(\alpha_h)}{R_\ast }   \dfrac{\partial \tilde  \alpha }{\partial \xi} + \dfrac{\hat \mu_c }{R_\ast^2} \dfrac{\partial^2 \tilde  v_c^1}{\partial \xi^2}+\dfrac{(\lambda_c+ \mu_c )}{ R_\ast } i\kappa \dfrac{\partial \tilde  v^2_c}{\partial \xi}   - \mu_c \kappa^2  \tilde  v_c^1
\\[.3cm]
0 & =  \Big( - \Sigma_c(\alpha_h)-\alpha_h \Sigma_c'(\alpha_h) + \lambda_c  \lambda_2  \Big)i\kappa \tilde \alpha
+ \dfrac{\mu_c \alpha_h}{R_\ast^2} \dfrac{\partial^2 \tilde  v_c^2}{\partial \xi^2}
\\
& + (\mu_c+\lambda_c)\dfrac{\alpha_h i\kappa}{R_\ast}\dfrac{\partial \tilde v_c^1}{\partial \xi}
- \alpha_h \kappa^2 \hat \mu_c  \tilde v_c^2
\\[.3cm]
0 &= \dfrac{\alpha_h }{R_\ast^2} \dfrac{\partial}{\partial \xi}( \tilde v_w^2 -  \tilde v_c^2) - \dfrac{i\kappa \alpha_h}{R_\ast} (\tilde v_w^1 - \tilde v_c^1)  +\dfrac{i\kappa \lambda_2 \xi}{1-\alpha_h}\tilde \alpha
\end{cases}
\end{align}
for $(\xi,t)\in (0,1) \times (0,\infty)$, with 
$$\dfrac{d \tilde  R}{d t}(t) = \lambda_2 \tilde R(t)+\tilde  v_c^1(1,t),\qquad t\in (0,\infty),$$
subject to the boundary conditions:
\begin{align}
\begin{cases}
\tilde v_c^1\hspace{4pt}=\tilde v_c^2=\tilde v_w^1=0, \qquad &\mbox{ at } \xi=0,
\\
0 \hspace{8pt}= - \Sigma_c'(\alpha_h) \tilde \alpha + \dfrac{\hat \mu_c }{R_\ast } \dfrac{\partial \tilde v_c^1}{\partial \xi}+\lambda_c  i\kappa \tilde v_c^2,\qquad &\mbox{ at } \xi= 1,
\\[.2cm]
0  \hspace{8pt}= i\kappa \Big( \tilde v_c^1 + 2 \lambda_2  \tilde R \Big) + \dfrac{1}{R_\ast }\dfrac{\partial \tilde v_c^2}{\partial \xi}, \qquad &\mbox{ at } \xi = 1,
\\[.2cm]
\tilde v_w^2  = \dfrac{i \kappa \lambda_2 R_\ast}{1-\alpha_h} \tilde R+\tilde v_c^2, \qquad &\mbox{ at } \xi = 1,
\end{cases}
\end{align}
for $t\in (0,\infty)$.  

\section{Asymptotic results}

Since this system cannot be solved analytically, we obtain asymptotic expressions for large $t$ to classify the stability of the one dimensional solution.  For that, we  restrict ourselves to the case where the base state solution is one of a growing tumour, i.e. $\lambda_2>0$ and therefore
$$\epsilon = \epsilon(t) :=: \dfrac{1}{R_\ast(t)} = \dfrac{1}{R_0 e^{\lambda_2 t}} \longrightarrow 0, \qquad \text{ when } \qquad t\rightarrow \infty.$$ 
Motivated by this limit, $\epsilon(t)$ will be used as a small parameter for $t$ large in the asymptotic approximations that follow. 

Once this asymptotic behaviour is characterised we can compare these to the one dimensional base state solution in order to, ultimately, reach a condition for instability of the growing solutions of the limit case with negligible nutrient uptake and cell drag.  Details of this derivation can be found in \cite{genovese2017asymptotic}.  The large-time limit is of singular perturbation type, necessitating the application of the method of matched asymptotic expansions.  

\subsection{Outer solution}

We first characterised the long time behaviour of the solution to the system \eqref{eq:limit2d} at a position $\xi\in [0,1]$ away from any possible boundary layers.  Notice that the first three equations of this system pertain only to $\{\tilde \alpha, \tilde v_c^1, \tilde v_c^2\}$.  Thus, we may decouple these equations to solve for $\{\tilde \alpha, \tilde v_c^1, \tilde v_c^2\}$ first and then for $\{\tilde v_w^1, \tilde v_w^2\}$.

Suppose that $\{ a_0, a_1, a_2 \}$ are decay rates of $\tilde \alpha$, $\tilde v_c^1$, and $\tilde v_c^2$, respectively, i.e.
$$\tilde \alpha(\xi,t) = e^{a_0t} \overline \alpha(\xi), \qquad 
\tilde v_c^1 (\xi,t) = e^{a_1t} \overline v_c^1 (\xi), \qquad 
\tilde v_c^2 (\xi,t) = e^{a_2t} \overline v_c^2 (\xi),$$
with 
$$\overline \alpha(\xi) =O(1), \qquad \overline v_c^1(\xi) =O(1), \qquad \overline v_c^2(\xi) =O(1), \qquad t\rightarrow 0,$$ 
for all relevant $\xi$ (away from any possible boundary layers).  By substituting these into the first three equations in \eqref{eq:limit2d}, we may prove via case-by-case logic, that there is only one set of values $\{a_0,a_1,a_2\}$ that balances these equations for large time and yields non-trivial solutions for $\{\overline \alpha, \overline v_c^1, \overline v_c^2\}$.  Then, using these results we applied the same reasoning to find $\{\tilde v_w^1, \tilde v_w^2\}$ using the last two equations of the system. 

This way, we obtained that the non-trivial outer solution for $\xi$ away from any possible boundary layers is, at first order, 
\begin{align}\label{eq:outersolutiondecay1}
\tilde \alpha (\xi,t) & = e^{(\gamma_0 -\lambda_2) t}  \overline \alpha (\xi),\\
\tilde v_c^1 (\xi,t) & = \dfrac{\gamma_1}{R_0 }   e^{(\gamma_0 -2\lambda_2) t} \dfrac{d \overline \alpha }{d \xi}(\xi),\label{eq:outersolutiondecay2}\\
\tilde v_c^2 (\xi,t) & = \gamma_3 e^{(\gamma_0-\lambda_2) t}  \overline \alpha(\xi),\label{eq:outersolutiondecay3}\\
\tilde v_w^1(\xi,t) & =  \dfrac{\lambda_2 R_0}{\alpha_h ( 1-\alpha_h)} e^{\gamma_0 t}\xi \overline \alpha (\xi),\label{eq:outersolutiondecay4}\\
\tilde v_w^2(\xi,t)  & = -  e^{(\gamma_0-\lambda_2) t} \Bigg [ \dfrac{\lambda_2}{i \kappa \alpha_h ( 1- \alpha_h)^2}\bigg(1+\xi \dfrac{d}{d\xi} \bigg) \overline \alpha(\xi) + \dfrac{\alpha_h \gamma_3 }{1-\alpha_h} \overline \alpha(\xi)  \Bigg ]\label{eq:outersolutiondecay5}
\end{align}
and it satisfies
\begin{align*}
\begin{cases}
\dfrac{\partial \tilde \alpha }{\partial t} &=  \Big[\dfrac{\partial S_c}{\partial \alpha}(\alpha_h, C_\infty) -\lambda_2 \Big]\tilde \alpha  -\alpha_h i\kappa \tilde v_c^2  \\[.3cm]
0 &= \dfrac{\lambda_2}{1-\alpha_h}\Big(1+\xi \dfrac{\partial }{\partial \xi}\Big) \tilde  \alpha  +\dfrac{ 1-\alpha_h}{R_\ast} \dfrac{\partial \tilde  v_w^1 }{\partial \xi } +\alpha_h i\kappa \tilde  v_c^2 + (1-\alpha_h) i\kappa  \tilde  v_w^2\\[.3cm]
0 & = - \Sigma_c'(\alpha_h)   \dfrac{\partial \tilde  \alpha }{\partial \xi} +(\lambda_c+ \mu_c ) i\kappa \dfrac{\partial \tilde  v^2_c}{\partial \xi}   - \mu_c \kappa^2 R_\ast   \tilde  v_c^1 \\[.3cm]
0 & =  \Big( - \Sigma_c(\alpha_h)-\alpha_h \Sigma_c'(\alpha_h) + \lambda_c  \lambda_2  \Big)i\kappa\tilde \alpha - \alpha_h \kappa^2 \hat \mu_c  \tilde v_c^2. \\[.3cm]
0 &=  - \dfrac{i\kappa \alpha_h}{R_\ast} \tilde v_w^1 +\dfrac{i\kappa \lambda_2 \xi}{1-\alpha_h}\tilde \alpha,
\end{cases}
\end{align*}
where
$$
\begin{array}{c}
\gamma_0  = \dfrac{\partial S_c}{\partial \alpha}(\alpha_h,C_\infty) - \dfrac{\alpha_h\Sigma_c'(\alpha_h)}{\hat \mu_c} +\lambda_2 \Big(\dfrac{\lambda_c}{\hat \mu_c} -1\Big),\\[.5cm]
\gamma_1  = -\dfrac{\Sigma_c'(\alpha_h)}{\mu_c \kappa^2}-\dfrac{ (\lambda_c+\mu_c)( - \Sigma_c(\alpha_h)-\alpha_h \Sigma_c'(\alpha_h) + \lambda_c  \lambda_2  )}{\mu_c \kappa^2 \hat \mu_c \alpha_h},\\[.5cm]
\gamma_3  = -  \dfrac{( - \Sigma_c(\alpha_h)-\alpha_h \Sigma_c'(\alpha_h) + \lambda_c  \lambda_2 )}{\hat \mu_c \alpha_h i \kappa}.
\end{array}
$$
Notice that $\gamma_0$, $\lambda_2$, and therefore all the exponential decay rates specified for the outer solution, interestingly enough, do not depend on the wavelength $\kappa$ of the perturbation in question. 

\subsection{Boundary layers}

\begin{figure}[H]
\centering
\begin{minipage}{.49\textwidth}
\centering
\includegraphics[width=\textwidth]{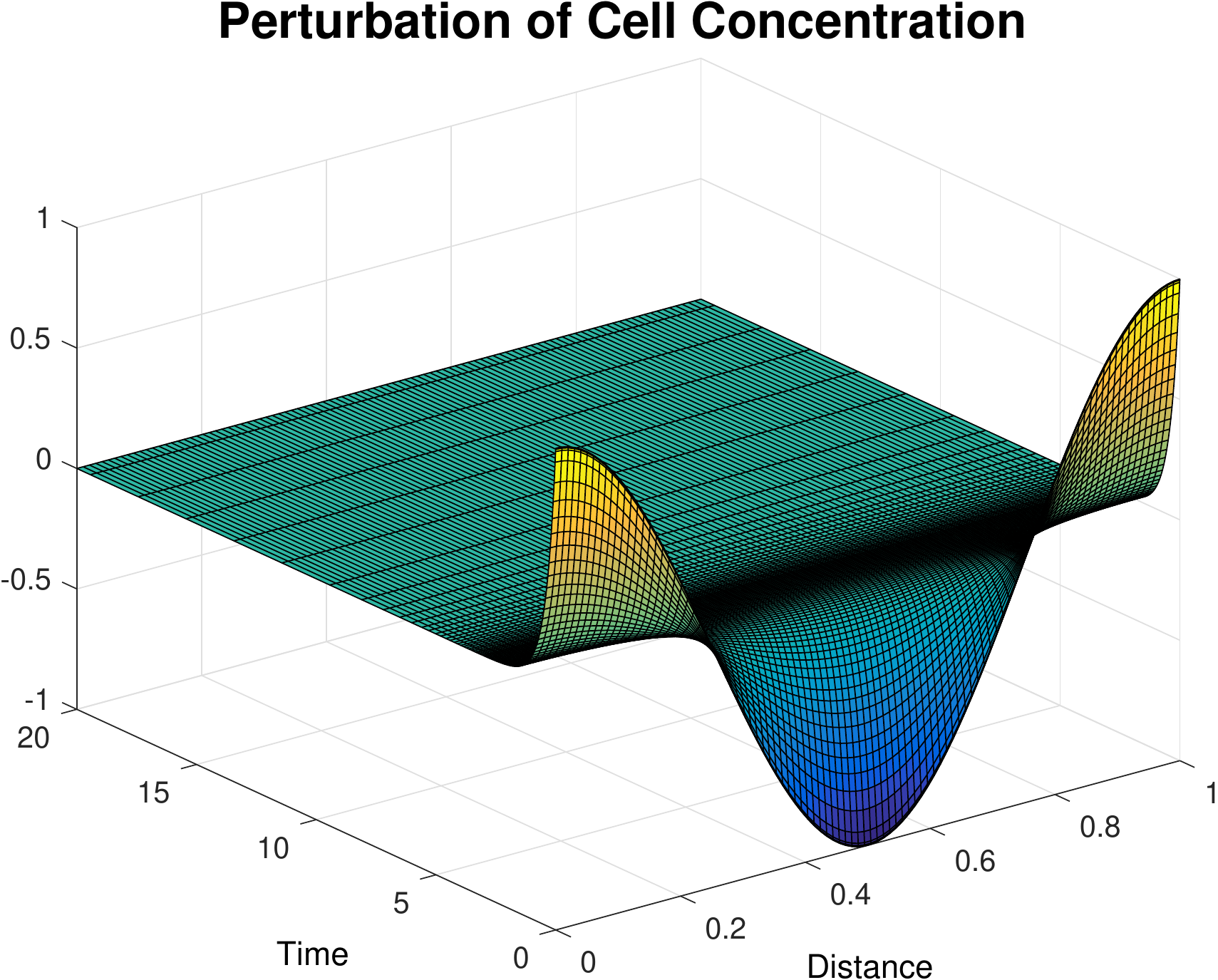}
\end{minipage}
\hfill
\begin{minipage}{.49\textwidth}
\centering
\includegraphics[width=\textwidth]{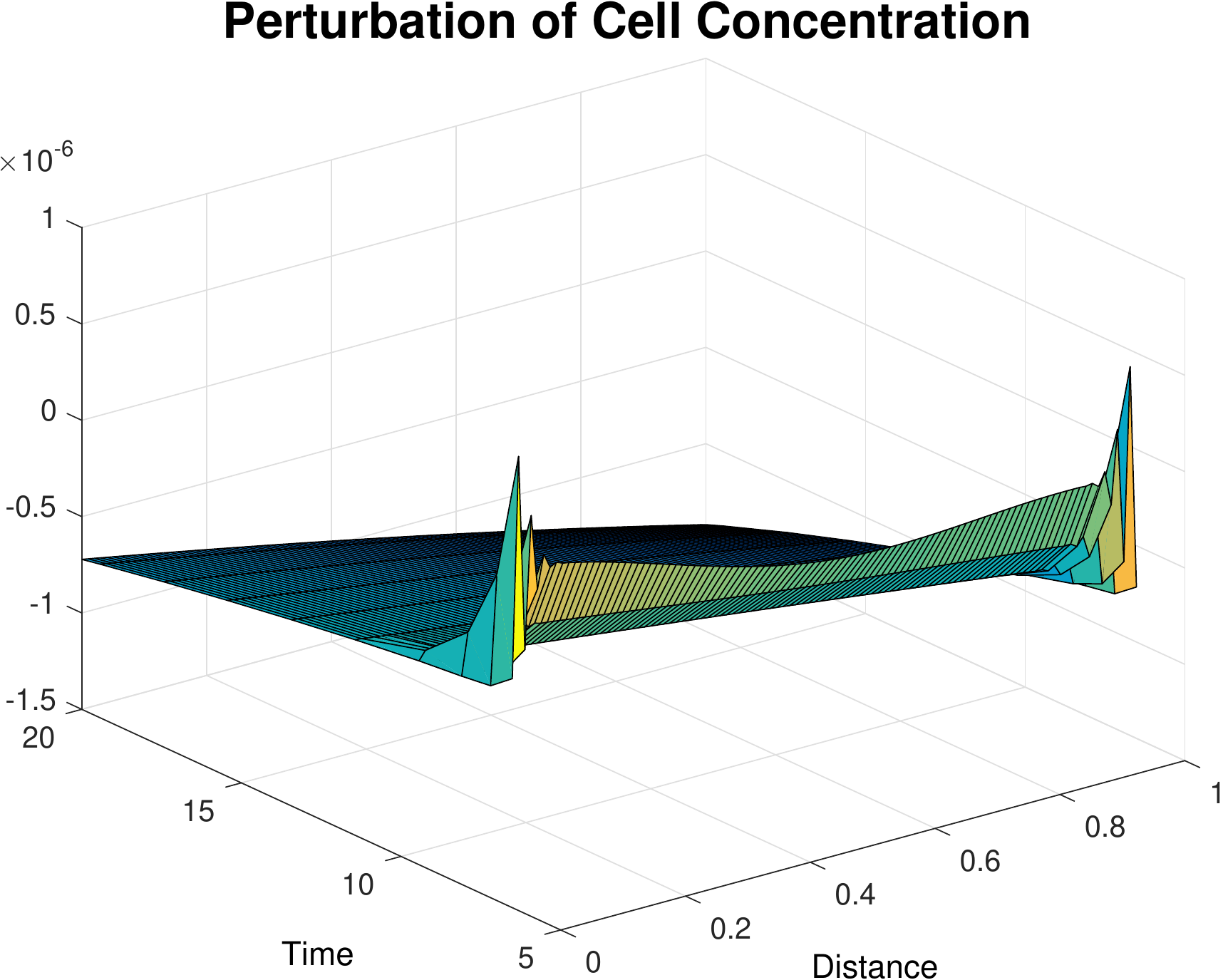}
\end{minipage}
\caption{Perturbation of the cell concentration $\tilde \alpha(\xi,t)$ for $t\in [0,20]$ (on the left) and $t\in[5,20]$ (on the right) showing rapid change near inner and outer boundaries $\xi=0$ and $\xi=1$, indicating the position of boundary layers.}
\label{fig:1}
\end{figure}

Now, in order to find the location of possible boundary layers, we must be able to find if there are any regions where rapid change may be observed.  For that, we tested different time frames of our solution numerically and discovered the following behaviour in $\tilde \alpha$ (see figure \ref{fig:1}).  

The graph on the right represents the same data as the figure on the left, except that the figure on the right shows the behaviour for $\tilde \alpha$ for $t>5$ rather than $t>0$, as seen on the left.  If we focus our attention on the graph on the right, we see that there are regions of rapid change near $\xi=0$ and $\xi=1$.  This, in turn, suggests that there is a boundary layer near each of these boundaries.

\subsection{Inner solution near free boundary}

Recall that, in the previous section, we have seen numerical results that suggest that there might be boundary layers near boundaries at $\xi=0$ and $\xi=1$, which we will call the inner and outer boundary, respectively.  Having said that, in this section, we will analytically characterise an inner solution near the outer boundary. Indeed, by considering all possible asymptotic balances for our equations near $\xi=1$, we obtain four possible non-trivial balances, which may be narrowed down to only one possibility that is non-trivial and allows for matching to be performed between the outer solution and the inner solution.  

For this, consider the following change of variable
$$\xi =1 - \dfrac{X}{R_\ast(t)^\beta}, \qquad \beta>0, \qquad X=O(1)\quad\text{ as }\quad t\rightarrow \infty.$$
This will allow us to focus within the boundary layer, a thin strip near the outer boundary $\xi=1$.  By the product rule, we then have 
\begin{align*}
\dfrac{\partial}{\partial t} & \longrightarrow \dfrac{\partial}{\partial t} + \beta \lambda_2 X \dfrac{\partial}{\partial X} \\[.1cm]
\dfrac{\partial}{\partial \xi} & \longrightarrow -R_\ast(t)^\beta \dfrac{\partial}{\partial X}\\[.1cm]
\dfrac{\partial^2}{\partial \xi^2} & \longrightarrow R_\ast(t)^{2\beta} \dfrac{\partial^2}{\partial X^2}.
\end{align*}

Recall also that the outer solutions all had a common coefficient: $e^{\gamma_0 t}$.   Thus, to facilitate matching between the outer and inner solutions, we will choose to assume that our inner solution is of the form 
$$\tilde \alpha(X,t) = e^{(\gamma_0 + a_0) t} A(X),\quad \tilde v_c^1(X,t) = e^{(\gamma_0 +a_1) t} V_c^1(X),\quad \tilde v_c^2(X,t) = e^{(\gamma_0+ a_2) t} V_c^2(X),$$
$$\tilde v_w^1(X,t) = e^{(\gamma_0+b_1) t} V_w^1(X), \quad \tilde v_w^2(X,t) = e^{(\gamma_0 +b_2 )t} V_w^2(X).$$
By substituting these into the system above and exploring the exponential pre-factors of each of its terms and their behaviour when $t\rightarrow \infty$, it is possible to prove through a case-by-case logic (similar to the one used before) that the only possible non-trivial asymptotic balances happen when
$$a_0= a_1 =a_2, \qquad b_1-a_0 - \lambda_2 = b_2 - a_0 - \lambda_2=0,$$
where the asymptotic balance is given by
\begin{empheq}[left=\empheqlbrace]{align}
0&=  \Big[\dfrac{\partial S_c}{\partial \alpha}(\alpha_h, C_\infty)-\gamma_0 -a_0 -\lambda_2 \Big( 1+ X\dfrac{d}{dX} \Big)\Big]A \nonumber \\
&+ \alpha_h \dfrac{ d V^1_{c}}{d X} -\alpha_h i\kappa V_c^2  \label{eq:AAI4_1} \\[.1cm] 
0 &= - \dfrac{\lambda_2 R_0 }{(1-\alpha_h)^2}\dfrac{d A}{d  X}  - \dfrac{d  V_w^1 }{d  X } + i\kappa V_w^2 \label{eq:AAI4_4}\\[.1cm] 
0 & = \Sigma_c'(\alpha_h) \dfrac{d A}{d X} + \hat \mu_c \dfrac{d^2 V_c^1}{d X^2}- (\lambda_c+ \mu_c ) i\kappa \dfrac{d V^2_c}{d X}   - \mu_c \kappa^2  V_c^1 \label{eq:AAI4_2}\\[.1cm] 
0 & =  \Big( - \Sigma_c(\alpha_h)-\alpha_h \Sigma_c'(\alpha_h) + \lambda_c  \lambda_2  \Big)  i\kappa A
+  \mu_c \alpha_h  \dfrac{d^2 V_c^2}{d X^2} \nonumber \\
  & - (\mu_c+\lambda_c) \alpha_h i\kappa \dfrac{d V_c^1}{d X}
- \alpha_h \kappa^2 \hat \mu_c  V_c^2  \label{eq:AAI4_3}\\[.1cm] 
0 &= -\dfrac{d V_w^2}{d  X}   - i\kappa  V_w^1 +\dfrac{i\kappa\lambda_2 R_0 }{\alpha_h(1-\alpha_h)} A, \label{eq:AAI4_5}
\end{empheq}
with conditions at $X=0$:
\begin{align*}
\begin{cases}
0  &= - \Sigma_c'(\alpha_h) A - \hat \mu_c \dfrac{d V_c^1}{d X}+\lambda_c  i\kappa V_c^2,\\[.3cm] 
0 &= \dfrac{(\lambda_2-\gamma_0-a_0)}{i\kappa } \dfrac{d V_c^2}{d X}+ (\lambda_2 +\gamma_0+a_0)V_c^1\\[.3cm]
V_w^2 & =\dfrac{R_0}{2(1-\alpha_h)}\Big[ \dfrac{d V_c^2}{dX}- i\kappa V_c^1\Big].  
\end{cases}
\end{align*}

Consider initially the first three equations \eqref{eq:AAI4_1}-\eqref{eq:AAI4_3} and the corresponding set $\{ A, V_c^1, V_c^2\}$.  Notice  that the asymptotic balance of these equations preserves all the terms in the original system.  Therefore, it is no surprise that we cannot solve this system by using basic ordinary differential equation methods as we have done in previous cases.  However, we can begin to tackle equations \eqref{eq:AAI4_1}-\eqref{eq:AAI4_3} by defining $\{F_1, F_2\}$ as  
$$F_1 = \dfrac{d V_c^1}{dX} -i\kappa V_c^2, \qquad F_2 =\dfrac{d V_c^2}{dX} +i\kappa V_c^1.$$
By substituting these into \eqref{eq:AAI4_1} and \eqref{eq:AAI4_2}, we can find expressions for $\{F_1, F_2\}$ in terms of $A$ which we can then substitute into \eqref{eq:AAI4_3} to obtain a third order ordinary differential equation for $A$:
\begin{equation}\label{eq:AAI4_7}
0 =  \lambda_2 X\dfrac{d^3 A}{dX^3}+\Big[a_0+\lambda_2(\dfrac{\lambda_c}{\hat \mu_c}+2)\Big]\dfrac{d^2 A}{dX^2}-\lambda_2 \kappa^2 X \dfrac{dA}{dX} - \kappa^2(a_0+\lambda_2)A.
\end{equation}
Then, through some algebraic manipulations, the expressions for $\{F_1, F_2\}$ in terms of $A$ allow us to recover $\{V_c^1,V_c^2\}$ by solving:  
\begin{align}
\begin{cases}\label{eq:AAI4_8}
\dfrac{d^2 V_c^1}{d X^2}-\kappa^2 V_c^1 & = H_1\\[.1cm]
V_c^2 -\dfrac{1}{i\kappa} \dfrac{dV_c^1}{dX} & = H_2,
\end{cases}
\end{align}
where $H_1$ and $H_2$ depend only on $A$:
\begin{align*}
H_1& =\dfrac{\mu_c + \lambda_c}{\mu_c \alpha_h} \bigg[ \Big(\dfrac{\partial S_c}{\partial \alpha}(\alpha_h, C_\infty)-\gamma_0-a_0 -2 \lambda_2 -\dfrac{\alpha_h \Sigma_c'(\alpha_h)}{\mu_c + \lambda_c}  \Big) \dfrac{dA}{dX}  - \lambda_2 X \dfrac{d^2 A}{d X^2}  \bigg]  \\[.2cm]
H_2& = \dfrac{1}{\alpha_h i\kappa} \Big(\dfrac{\partial S_c}{\partial \alpha}(\alpha_h, C_\infty)-\gamma_0 -a_0-\lambda_2  \Big)A - \dfrac{\lambda_2}{\alpha_h i \kappa } X \dfrac{dA}{dX}.
\end{align*}

However, this whole simplification process relies on solving $A$ from a third order ordinary differential equation, namely \eqref{eq:AAI4_7}, which cannot be easily solved analytically.  Therefore, we will approach the problem differently and find the asymptotic behaviour of $\{A,V_c^1,V_c^2\}$ when $X\rightarrow \infty$ in order to be able to perform matching between the outer and inner solutions to obtain a value for $a_0$ and consequently for $\{a_1,a_2,b_1,b_2\}$.  

For this, consider $\hat X = \epsilon X$, with $0<\epsilon\ll 1$.  Thus, equations \eqref{eq:AAI4_7} and \eqref{eq:AAI4_8} above may be found to be, at first order: 
$$0=\lambda_2 \hat X \dfrac{dA}{d \hat X}(\hat X) + (a_0+\lambda_2)A(\hat X),\quad
V_c^1(\hat X) = 0,$$
$$V_c^2(\hat X)   = \dfrac{1}{\alpha_h i\kappa} \Big[\dfrac{\partial S_c}{\partial \alpha}(\alpha_h, C_\infty)-\gamma_0 +\lambda_2 \Big(1-\hat X \dfrac{d }{d \hat X}\Big)\Big] A(\hat X).$$
Through traditional ordinary differential equation methods, these can be solved to find that 
$$A(\hat X) = C_1 \hat X^{-\dfrac{a_0+\lambda_2}{\lambda_2}}, \qquad V_c^1(\hat X) = 0, \qquad   V_c^2(\hat X)= \gamma_3 C_1 \hat X^{-\dfrac{a_0+\lambda_2}{\lambda_2}}.$$ 
Thus, as $X\rightarrow +\infty$,
$$A(X) \sim C_1 X^{-\dfrac{a_0+\lambda_2}{\lambda_2}}, \qquad V_c^1(X) \rightarrow 0, \qquad V_c^2(X) \sim \gamma_3 C_1 X^{-\dfrac{a_0+\lambda_2}{\lambda_2}}.$$
Now, by using this asymptotic behaviour and recasting it in terms of $(\xi,t)$ using $X=(1-\xi)R_\ast(t)$, we see that as $X\rightarrow +\infty$,
\begin{align*}
e^{a_0 t} A(X) & \sim e^{a_0 t} C_1 X^{-\dfrac{a_0+\lambda_2}{\lambda_2}} 
= C_1 [R_0(1-\xi)]^{-\dfrac{a_0+\lambda_2}{\lambda_2}} e^{-\lambda_2 t},\\
e^{a_0 t} V_c^1(X) & \rightarrow 0,\\
e^{a_0 t} V_c^2(X) & \sim e^{a_0 t} C_1 \gamma_3 X^{-\dfrac{a_0+\lambda_2}{\lambda_2}} 
 = \gamma_3 C_1 [R_0(1-\xi)]^{-\dfrac{a_0+\lambda_2}{\lambda_2}} e^{-\lambda_2 t}.
\end{align*}
By matching these to the outer solutions, we see that $a_0 = -2\lambda_2$ and
$$\overline \alpha(\xi) \sim C_1 R_0(1-\xi), \qquad \text{ and } \qquad \dfrac{d \overline \alpha}{d \xi}(\xi) \rightarrow 0$$
when $\xi \rightarrow 1$.  Furthermore, 
we see that 
%
$$a_0=a_1=a_2=-2\lambda_2, \qquad b_1=b_2 = a_0+\lambda_2 =-\lambda_2,$$
and thus we have characterised the decay rates of the inner solution near the outer boundary.  

If we apply the same process used above for $\{A,V_c^1,V_c^2\}$ now to $\{V_w^1,V_w^2\}$, we see that equations \eqref{eq:AAI4_4} and \eqref{eq:AAI4_5} imply that, at first order,   
$$V_w^1(X) \sim \dfrac{\lambda_2 R_0}{\alpha_h (1-\alpha_h)} C_1 X \qquad \text{ and }\qquad V_w^2(X) \rightarrow 0,$$
when $X\rightarrow \infty$ and matching can be performed due to values of $b_1$ and $b_2$ and the conditions previously derived on $\overline \alpha(\xi)$, namely
$$\overline \alpha(\xi) \sim C_1 R_0(1-\xi), \qquad \text{ and } \qquad \dfrac{d \overline \alpha}{d \xi}(\xi) \rightarrow 0,$$
when $\xi \rightarrow 1$.  

Now, notice that the two conditions above together imply that $C_1=0$ (via L'Hopital's Rule) and therefore that the asymptotic behaviour obtained here when $X\rightarrow \infty$ is actually trivial.  In order to remedy this, we must revisit equation \eqref{eq:AAI4_7}, now with $a_0=-2\lambda_2$, namely
$$0 =  X\dfrac{d^3 A}{dX^3}+\dfrac{\lambda_c}{\hat \mu_c}\dfrac{d^2 A}{dX^2}-\kappa^2 X \dfrac{dA}{dX} + \kappa^2 A.$$
Notice that this is a third order ordinary differential equation, so we expect three linearly independent solutions and not just one, as we obtained previously with the first order approximation.  Therefore, we will extract the other two solutions by applying a WKBJ expansion to $A(X)$.  

Initially, consider $\hat X = \epsilon X,$ with $0<\epsilon\ll 1.$  Then, the equation above becomes: 
$$0 = \epsilon^2 \Big (\hat X \dfrac{d^3 A}{d \hat X^3} + \dfrac{\lambda_c}{\hat \mu_c}\dfrac{d^2 A}{d \hat X^2}\Big) -\kappa^2 \Big( \hat X \dfrac{dA}{d \hat X}-A\Big).$$
Note that, at first order, we can retrieve the solution obtained in the previous section:  $A(\hat X)= C_1 \hat X$. Now, as usual in a WKBJ approximation, let 
$$A(\hat X, \epsilon ) = \exp\Big\{\dfrac{G(\hat X,\epsilon)}{\epsilon}\Big\},\quad \text{ where } \quad G(\hat X,\epsilon)= G_0(\hat X) + \epsilon G_1(\hat X) + O( \epsilon^2).$$ 
Through substitution into the equation above, we obtain 
\begin{align*}
0 &= \hat X  \Bigg[  \Bigg(\dfrac{d G_0}{d \hat X} \Bigg)^3 - \kappa^2 \dfrac{d G_0}{d \hat X} \Bigg]+ \epsilon \Bigg[\kappa^2 + \dfrac{\lambda_c}{\hat \mu_c}\Bigg(\dfrac{d G_0}{d \hat X} \Bigg)^2 + 3\hat X \dfrac{d G_0}{d \hat X} \dfrac{d^2 G_0}{d \hat X^2} \\
& + 3 \hat X \Bigg(\dfrac{d G_0}{d \hat X} \Bigg)^2 \dfrac{d G_1}{d \hat X} - \hat X \kappa^2 \dfrac{d G_1}{d \hat X}  \Bigg]+O(\epsilon^2).
\end{align*}
Equating terms in $O(1)$ and then $O(\epsilon)$, we see that either
$$\{G_0(\hat X) = \tilde C_0, \quad G_1(\hat X) = \ln(\hat X) + \tilde C_1\},$$
or
$$\{G_0(\hat X) = \pm \kappa \hat X + \tilde C_0, \quad G_1(\hat X) = \omega \ln(\hat X) + \tilde C_1\},$$
where $\omega = -\dfrac{1}{2}\big(1+\dfrac{\lambda_c}{\hat \mu_c}\big)$. Therefore, grouping the constants, we have 
$$A(\hat X, \epsilon ) = \exp \Big\{ \frac{G(\hat X,\epsilon)}{\epsilon} \Big\} = \exp \Big\{\frac{G_0(\hat X)}{\epsilon}+G_1(\hat X) \Big\} = C_1 \hat X,$$ 
or 
$$A(\hat X, \epsilon ) = \exp \Big\{ \frac{G(\hat X,\epsilon)}{\epsilon} \Big\} = \exp \Big\{\frac{G_0(\hat X)}{\epsilon}+G_1(\hat X) \Big\}  = C_0 \hat X^{\omega} \exp \Big\{\dfrac{\pm \kappa\hat X}{\epsilon}\Big\}.$$
Notice that the first solution is the same as found previously, by using a first order approximation.  Thus, the latter equation corresponds to the other two solutions for which we searched.  Therefore, since the equation for $A$ is a third order ordinary differential equation, we thus conclude that $A(X)$ satisfies
\begin{equation}\label{eq:AAI4_9}
A(X) \sim C_1 X + C_2 X^{\omega}e^{\kappa X} + C_3 X^{\omega}e^{-\kappa X}, \qquad \omega = -\dfrac{1}{2}\big(1+\dfrac{\lambda_c}{\hat \mu_c}\big)
\end{equation} 
when $X\rightarrow \infty$ for some constants $C_1, C_2, C_3$.  However, for matching to be possible, we must have $C_2=0$.  

Recall now the equations in \eqref{eq:AAI4_8} for $\{V_c^1, V_c^2\}$ were obtained via algebraic manipulation and the definition of two auxiliary functions $\{ F_1, F_2\}$ that depended on $\{V_c^1, V_c^2\}$.  Similarly, we take an analogous approach to recast equations \eqref{eq:AAI4_4} and \eqref{eq:AAI4_5} to find $\{V_w^1,V_w^2\}$.  To be more specific, we define
$$F_3 = \dfrac{d V_w^1}{dX} -i\kappa V_w^2, \qquad F_4 =\dfrac{d V_w^2}{dX} +i\kappa V_w^1,$$
and, through some algebraic manipulation, recast equations \eqref{eq:AAI4_4} and \eqref{eq:AAI4_5} as
\begin{align}
\begin{cases}\label{eq:AAI4_10}
\dfrac{d^2 V_w^2}{d X^2}-\kappa^2 V_w^2 & = \dfrac{i \kappa R_0 \lambda_2}{\alpha_h (1-\alpha_h)^2} \dfrac{dA}{dX}\\[.1cm]
V_w^1 +\dfrac{1}{i\kappa} \dfrac{d V_w^2}{dX} & = \dfrac{\lambda_2 R_0}{\alpha_h (1-\alpha_h)}A .     
\end{cases}
\end{align}

Therefore, using the expression for $A$ in \eqref{eq:AAI4_9}, we can find the asymptotic behaviour of $\{V_c^1, V_c^2, V_w^1, V_w^2\}$ via equations \eqref{eq:AAI4_8} and \eqref{eq:AAI4_10}.  For this, we must solve two second order ordinary differential equations for $\{V_c^1,V_w^2\}$ and obtain $\{V_c^2,V_w^1\}$ by substitution of $\{A,V_c^1,V_w^2\}$ into the remaining equations.  

Note also that the solution to the homogeneous version of the ordinary differential equations for $\{V_c^1,V_w^2\}$ will have both decaying and growing exponential terms of the type $e^{-\kappa X}$ and $e^{\kappa X}$, respectively.  However, for matching to be possible, we cannot allow the existence of the growing exponential term of the type $e^{\kappa X}$.  If this is taken into consideration, we can simplify our calculations by assuming any pre-factor of $e^{\kappa X}$ is null.  Thus, by inserting the expression for $A(X)$ into the equations above, we obtained that
\begin{align*}
A(X) & \sim  C_1 X +C_2 X^{\omega} e^{\kappa X}+C_3 X^{\omega} e^{-\kappa X}\\
& \Longrightarrow C_2 = 0\\
V_c^1(X) & \sim -C_1 \gamma_1 + C_4 e^{\kappa X}+C_5 e^{-\kappa X}+C_3 X^{\omega} e^{-\kappa X} p_1^2(X)\\
& \Longrightarrow C_4 = 0\\
V_c^2(X) & \sim \gamma_3 C_1 X + C_5 i e^{-\kappa X} +C_3 X^{\omega-1} e^{-\kappa X} p_2^3(X)\\
V_w^1(X) & \sim \dfrac{\lambda_2 R_0}{\alpha_h ( 1-\alpha_h)} C_1 X + C_6 e^{-\kappa X}+C_3 X^{\omega-1} e^{-\kappa X} p_3^2(X) \\[.1cm]
V_w^2(X) & \sim \dfrac{\lambda_2 R_0 C_1}{i \kappa \alpha_h ( 1-\alpha_h)^2}  +C_6 i e^{-\kappa X}+C_7 e^{\kappa X}+C_3 X^{\omega} e^{-\kappa X}p_4^1(X)\\
& \Longrightarrow C_7 = 0,
\end{align*}
when $X\rightarrow \infty$.  Here we used the following notation:
$$p_i^n (X) = A_{i,0}+A_{i,1}X+A_{i,2}X^2 ... + A_{i,n}X^n, \qquad i=1,2,3,4, \qquad n\in \mathbf{N},$$
where every $A_{i,j}$ is a constant that may be found explicitly in terms of the parameters of our system.  Because these polynomials appear in terms that decay exponentially when $X\rightarrow \infty$, there is no need to make the expressions for $A_{i,j}$ explicit here as these are negligible terms for the purpose of matching. 

In the process used to derive the expressions above, we also used the first order term of the following asymptotic expansion: 
\begin{align*}
\int^X e^{-2\kappa s} s^\eta ds & \sim - \dfrac{X^\eta  e^{-2\kappa X}}{2\kappa} \Big( 1+ \dfrac{\eta }{2\kappa}\dfrac{1}{X}+\dfrac{\eta (\eta -1)}{(2\kappa)^2}\dfrac{1}{X^2}\\
 & +\dfrac{\eta (\eta -1)(\eta -2)}{(2\kappa)^3}\dfrac{1}{X^3} + O\Big(\Big(\dfrac{1}{X}\Big)^4\Big) \Big),
\end{align*}
when $X\rightarrow \infty$, which can be proven by integrating by parts repeatedly (See \cite{temme2014asymptotic}). 

The terms that grow exponentially must be discarded in order to allow us to perform matching when $X\rightarrow \infty$, and therefore $C_2=C_4=C_7=0$.  This way, no exponentially growing term is carried along in the derivation of $V_c^1$, $V_c^2$, $V_w^1$, and $V_w^2$ and we may obtain $3$ boundary conditions for the inner solution when $X\rightarrow \infty$.  Since there are $3$ boundary conditions already at $X=0$, this implies that we have a total of $6$ boundary conditions that may be applied to the inner solution at $X=0$ and when $X\rightarrow \infty$ in order to perform numerical simulations.  

Notice that thus we have $4$ degrees of freedom, corresponding to the constants $\{ C_1,C_3,C_5,C_6 \}$.  However, $\{ C_3,C_5,C_6 \}$ do not play a role in matching, since these are the coefficients of terms that decay exponentially when $X\rightarrow \infty$. 

That being said, we will now find a condition for $C_1$ by matching the inner solution to the outer solution when $X\rightarrow \infty$.  Indeed, from the expressions above, we see that when $X\rightarrow \infty$, at first order
\begin{align*}
A(X) & \sim  C_1 X \\
V_c^1(X) & \sim - C_1 \gamma_1 \\
V_c^2(X) & \sim C_1 \gamma_3 X \\
V_w^1(X) & \sim \dfrac{\lambda_2 R_0}{\alpha_h ( 1-\alpha_h)} C_1 X \\[.1cm]
V_w^2(X) & \sim \dfrac{\lambda_2 R_0 C_1}{i\kappa \alpha_h ( 1-\alpha_h)^2}.
\end{align*}
Notice that the asymptotic far-field behaviour found for $A$, $V_c^2$, and $V_w^1$ is the same as found previously, by using a first order approximation.  However, instead of having 
$$V_c^1(X) \rightarrow 0, \qquad \text{ and } \qquad V_w^2(X) \rightarrow 0$$
when $X\rightarrow \infty$, we now have 
$$V_c^1(X) \sim - C_1 \gamma_1, \qquad \text{ and } \qquad V_w^2(X) \sim \dfrac{\lambda_2 R_0 C_1}{i\kappa \alpha_h ( 1-\alpha_h)^2}$$
when $X\rightarrow \infty.$  If we match these to the outer solution when $\xi\rightarrow 1$, we obtain that for matching to be possible we must have 
$$C_1 = -\dfrac{1}{R_0} \lim_{\xi \rightarrow 1} \dfrac{d \overline \alpha}{d \xi}(\xi), \qquad \text{ and } \qquad \overline \alpha(\xi) \sim C_1 R_0 (1-\xi),$$
when $\xi\rightarrow 1$, where $\overline \alpha(\xi)$ is part of the outer solution.  

As in the previous section, notice again that all the exponential decay rates specified for the inner solution near the outer boundary, interestingly enough, do not depend on the wavelength $\kappa$ of the perturbation in question. 

\subsection{Inner solution near fixed boundary}

Finally, we characterised the inner solution in the boundary layer near the inner boundary $\xi=0$.  For that, we have performed a process very similar to the one applied to the characterisation of the inner solution near the outer boundary.  Most changes stem from different rescaling on $\xi$, namely 
$$\xi =\dfrac{x}{R_\ast(t)^\beta}, \quad x=O(1)\text{ as }t\rightarrow \infty, \quad \beta>0.$$
Indeed, we now want to consider a thin strip near $\xi=0$, rather than $\xi=1$, which justifies the rescaling above. 

Following the same process as in the previous section, we found that the inner solution near the inner boundary $\xi=0$ satisfies at order one:  
\begin{align*}
\tilde \alpha (x,t) & = e^{(\gamma_0 -2 \lambda_2) t}  A(x)\\
\tilde v_c^1 (x,t) & = e^{(\gamma_0 -2 \lambda_2) t}  V_c^1(x)\\
\tilde v_c^2 (x,t) & = e^{(\gamma_0 -2 \lambda_2) t}  V_c^2(x)\\
\tilde v_w^1 (x,t) & = e^{(\gamma_0 - 2\lambda_2) t}  V_w^1(x)\\
\tilde v_w^2 (x,t) & = e^{(\gamma_0 - 2\lambda_2) t}  V_w^2(x)
\end{align*}
where 
$$\xi=\dfrac{x}{R_\ast(t)}$$
and $A(x),$ $V_c^1(x)$, $V_c^2(x)$, $V_w^1(x)$, and $V_w^2(x)$ satisfy
\begin{align*}
\begin{cases}
0&=  \Big[\dfrac{\partial S_c}{\partial \alpha}(\alpha_h, C_\infty)-\gamma_0 +\lambda_2 \Big( 1- x\dfrac{d}{dx} \Big)\Big]A- \alpha_h \dfrac{ d V^1_{c}}{d x} -\alpha_h i \kappa  V_c^2  \\[.2cm] 
0 & = -\Sigma_c'(\alpha_h) \dfrac{d A}{d x} + \hat \mu_c \dfrac{d^2 V_c^1}{d x^2}+ (\lambda_c+ \mu_c ) i \kappa   \dfrac{d V^2_c}{d x}   - \mu_c \kappa^2  V_c^1 \\[.2cm] 
0 & =  \Big( - \Sigma_c(\alpha_h)-\alpha_h \Sigma_c'(\alpha_h) + \lambda_c  \lambda_2  \Big)  i \kappa  A
+  \mu_c \alpha_h  \dfrac{d^2 V_c^2}{d x^2}\\
& + (\mu_c+\lambda_c) \alpha_h i \kappa  \dfrac{d V_c^1}{d x}
- \alpha_h \kappa^2 \hat \mu_c  V_c^2  \\[.2cm] 
0 &=  \dfrac{\lambda_2}{1-\alpha_h} \big( 1+x \dfrac{d}{dx} \big) A  +\alpha_h \dfrac{d V_c^1 }{d x }+\alpha_h i \kappa V_c^2 \\
& + (1-\alpha_h) \dfrac{d  V_w^1 }{d  x } + (1-\alpha_h) i \kappa   V_w^2\\[.2cm] 
0 &= \dfrac{\lambda_2 i \kappa }{1-\alpha_h}x A +\alpha_h i \kappa  V_c^1-\alpha_h \dfrac{d V_c^2 }{d x }  - \alpha_h i \kappa   V_w^1+ \alpha_h \dfrac{d V_w^2}{d  x},
\end{cases}
\end{align*}
and $V_c^1(0)=V_c^2(0)=V_w^1(0)=0$.

Far field behaviour when $x\rightarrow \infty$ is, at first order
\begin{align*}
A(x) & \sim  D_1 x \\
V_c^1(x) & \sim D_1 \gamma_1 \\
V_c^2(x) & \sim D_1 \gamma_3 x \\
V_w^1(x) & \sim \dfrac{\lambda_2 }{\alpha_h ( 1-\alpha_h)} D_1 x^2\\
V_w^2(x) & \sim D_1 B_1 x,
\end{align*}
where matching yields
$$D_1 = \dfrac{1}{R_0} \lim_{\xi \rightarrow 0} \dfrac{d \overline \alpha}{d \xi}(\xi), \qquad \text{ and } \qquad \overline \alpha(\xi) \sim D_1 R_0 \xi,$$
when $\xi \rightarrow 0$.  Details of this derivation can be found in \cite{genovese2017asymptotic}.

As in the previous sections, notice again that all the exponential decay rates specified for the inner solution near the inner boundary, interestingly enough, do not depend on the wavelength $\kappa$ of the perturbation in question. 

At this point, notice that the asymptotic characterisation of the outer and inner solutions has been obtained. 

Additionally, the decays found in this characterisation were validated numerically by simulating the system satisfied by the perturbations and comparing decays to those characterised in the previous sections, for various regions of the tumour.  The interested reader may find more details in \cite{genovese2017asymptotic}.

\section{Summary}

The non-homogeneous base state solutions are unstable to two dimensional perturbations if one of the following has an infinite limit for some $(x,y)\in (0, R_\ast(t))\times \mathcal{R}$ when $t\rightarrow \infty$:
$$\Bigg|\dfrac{\hat \alpha(x,y,t)}{\alpha_h}\Bigg|, \qquad  \dfrac{ || \mathbf{\hat  v}_c(x,y,t)|| }{||\mathbf{v}_{cs}(x)||},\qquad  \dfrac{|| \mathbf{\hat v}_w(x,y,t)|| }{||\mathbf{v}_{ws}(x)||}, \qquad  \Bigg|\dfrac{\hat R(y,t)}{R_\ast(t)}\Bigg|,$$
%
where $||\cdot || $ is the euclidean norm.  Equivalently, since the base state solution is one dimensional, it is linearly unstable if one of the following is infinite for some $(x,y)\in (0, R_\ast(t))\times \mathcal{R}$  when $t\rightarrow \infty$:
$$\Bigg|\dfrac{\hat \alpha(x,y,t)}{\alpha_h}\Bigg|, \qquad  \Bigg|\dfrac{\hat v_c^j(x,y,t)}{v_{cs}(x)}\Bigg|, \qquad  \Bigg|\dfrac{\hat v_w^j(x,y,t)}{v_{ws}(x)}\Bigg|, \qquad \Bigg|\dfrac{\hat R(y,t)}{R_\ast(t)}\Bigg|, \qquad j=1,2.$$

By inserting the definition of the base state solutions into the ratios above and using the $(\xi,t)$, $(X,t)$ and $(x,t)$ formulations defined in the previous section, through some algebraic manipulation, it is thus possible to prove that the base state solution is linearly unstable to two dimensional perturbations if one of the following is infinite for some $\xi\in (0,1)$, $X,x \in [0,\infty)$  when $t\rightarrow \infty$:
\begin{table}[H]
\centering
\begin{tabular}{p{2.5cm} p{3cm} p{3cm} }
$|\tilde \alpha(\xi,t)|$ & $|\tilde v_c^j(\xi,t)| e^{-\lambda_2 t}$ & $|\tilde v_w^j(\xi,t)|e^{-\lambda_2 t}$ \\ 
$|\tilde \alpha(X,t)|$ & $|\tilde v_c^j(X,t)| e^{-\lambda_2 t}$ & $|\tilde v_w^j(X,t)| e^{-\lambda_2 t}$ \\ 
 $|\tilde \alpha(x,t)|$ & $|\tilde v_c^j(x,t)|$ & $|\tilde v_w^j(x,t)|$ 
\end{tabular}
\end{table}
\noindent for $j=1,2$ or 
$$\lim_{t\rightarrow \infty} |\tilde R(t)| e^{-\lambda_2 t} = \infty.$$

By inserting the asymptotic characterisation of the perturbations obtained in the previous section into the expressions above, we can find a condition for instability of the growing base state solutions in terms of $\gamma_0$ and $\lambda_2$.  By doing this, we may conclude that for the base state solution to be unstable, one of the following exponentials must tend to $\infty$ when $t\rightarrow \infty$:
$$e^{(\gamma_0-\lambda_2)t}, \qquad e^{(\gamma_0-2\lambda_2)t}, \qquad e^{(\gamma_0-3\lambda_2)t}.$$
Thus, since $\lambda_2>0$, for instability we must have $\gamma_0-\lambda_2>0$.
\vspace{.5cm}

\subsection{Conclusion}

The growing base state solutions $(\alpha_h, \mathbf{v}_{cs}, \mathbf{v}_{ws}, \mathbf{R}_\ast)$ of the two dimensional limit case system with negligible nutrient uptake and cell drag obtained in section \ref{sec:3} are unstable to two dimensional perturbations when 
$$\gamma_0  -\lambda_2 = \dfrac{\partial S_c}{\partial \alpha}(\alpha_h,C_\infty) - \dfrac{\alpha_h\Sigma_c'(\alpha_h)}{\hat \mu_c} +\lambda_2 \Big(\dfrac{\lambda_c}{\hat \mu_c} -2\Big) >0,$$
where $\alpha_h$ and $\lambda_2$ are defined by
$$\hat\mu_c S_c(\alpha_h, C_\infty) = \alpha_h \Sigma_c(\alpha_h), \qquad \lambda_2 =\dfrac{\Sigma_c(\alpha_h)}{\hat\mu_c}.$$

\vspace*{6pt}

\noindent\textbf{Acknowledgements}\vspace*{3pt}

\noindent This work was produced with the financial support of the University of Nottingham, \textit{Coordena\c{c}\~ao de Aperfei\c{c}oamento de Pessoal de N\'ivel Superior} (CAPES), \textit{Funda\c{c}\~ao de Apoio a Pesquisa do Distrito Federal} (FAPDF), \textit{Conselho Nacional de Desenvolvimento Cient\'ifico e Tecnol\'ogico} (CNPq) and the University of Brasilia along with the collaboration of Prof. Daniele Avitabile from the Department of Mathematics, Vrije Universiteit Amsterdam, Netherlands.   

\end{document}

%% file: standard.tex
 \newtheoremstyle{theorem}{6pt}{6pt}{\rm}{}{\sffamily}{ }{ }{}
 \theoremstyle{theorem}

 \newtheoremstyle{algorithm}{6pt}{6pt}{\rm}{}{\sffamily}{ }{ }{}
 \theoremstyle{algorithm}

 \newtheoremstyle{lemma}{6pt}{6pt}{\rm}{}{\sffamily}{ }{ }{}
 \theoremstyle{lemma}

\newtheoremstyle{case}{6pt}{6pt}{\rm}{}{\sffamily}{. }{ }{}
 \theoremstyle{case}

 \newtheoremstyle{statement}{6pt}{6pt}{\rm}{}{\sffamily}{ }{ }{}
\theoremstyle{statement}

 \newtheoremstyle{corollary}{6pt}{6pt}{\rm}{}{\sffamily}{ }{ }{}
 \theoremstyle{corollary}

  \newtheoremstyle{definition}{6pt}{6pt}{\rm}{}{\sffamily}{ }{ }{}
 \theoremstyle{definition}

\newtheoremstyle{example}{6pt}{6pt}{\rm}{}{\sffamily}{ }{ }{}
\theoremstyle{example}

\newtheoremstyle{remark}{6pt}{6pt}{\rm}{}{\sffamily}{ }{ }{}
\theoremstyle{remark}

\newtheoremstyle{approximation}{6pt}{6pt}{\rm}{}{\sffamily}{ }{ }{}
\theoremstyle{approximation}

\newtheoremstyle{scheme}{6pt}{6pt}{\rm}{}{\sffamily}{ }{ }{}
\theoremstyle{scheme}

\newtheoremstyle{Algorithm}{6pt}{6pt}{\rm}{}{\sffamily}{ }{ }{}
\theoremstyle{Algorithm}

\newtheoremstyle{Assumption}{6pt}{6pt}{\rm}{}{\sffamily}{ }{ }{}
\theoremstyle{Assumption}

\newtheoremstyle{proposition}{6pt}{6pt}{\rm}{}{\sffamily}{ }{ }{}
\theoremstyle{proposition}

\newtheoremstyle{hypo}{6pt}{6pt}{\rm}{}{\sffamily}{ }{ }{}
 \theoremstyle{hypo}

  \newtheoremstyle{Step}{6pt}{6pt}{\rm}{}{}{ }{ }{}
 \theoremstyle{Step}